# Residual-driven online Generalized Multiscale Finite Element Methods

Eric T. Chung[*]    Yalchin Efendiev[†]    Wing Tat Leung[‡]

January 19, 2015

**Abstract**

The construction of local reduced-order models via multiscale basis functions has been an area of active research. In this paper, we propose online multiscale basis functions which are constructed using the offline space and the current residual. Online multiscale basis functions are constructed adaptively in some selected regions based on our error indicators. We derive an error estimator which shows that one needs to have an offline space with certain properties to guarantee that additional online multiscale basis function will decrease the error. This error decrease is independent of physical parameters, such as the contrast and multiple scales in the problem. The offline spaces are constructed using Generalized Multiscale Finite Element Methods (GMsFEM). We show that if one chooses a sufficient number of offline basis functions, one can guarantee that additional online multiscale basis functions will reduce the error independent of contrast. We note that the construction of online basis functions is motivated by the fact that the offline space construction does not take into account distant effects. Using the residual information, we can incorporate the distant information provided the offline approximation satisfies certain properties.

In the paper, theoretical and numerical results are presented. Our numerical results show that if the offline space is sufficiently large (in terms of the dimension) such that the coarse space contains all multiscale spectral basis functions that correspond to small eigenvalues, then the error reduction by adding online multiscale basis function is independent of the contrast. We discuss various ways computing online multiscale basis functions which include a use of small dimensional offline spaces.

## 1   Introduction

Solving real-world multiscale problems requires some type of model reduction due to disparity of scales. Many methods have been developed which can be classified as global [31, 36, 27] and local model reduction techniques [15, 38, 4, 6, 1, 16, 3, 30, 2, 34, 33, 29, 21, 23, 24, 27, 7, 12, 8, 9, 28]. Global model reduction techniques use global basis functions to construct reduced dimensional approximations for the solution space. These methods can involve costly offline constructions and lack local adaptivity. In this paper, our focus on the development of efficient local multiscale model reduction techniques that involve some local online computations.

[*]Department of Mathematics, The Chinese University of Hong Kong, Hong Kong SAR. This research is partially supported by the Hong Kong RGC General Research Fund (Project number: 400411).

[†]Department of Mathematics, Texas A&M University, College Station, TX; Numerical Porous Media SRI Center, King Abdullah University of Science and Technology (KAUST), Thuwal 23955-6900, Kingdom of Saudi Arabia

[‡]Department of Mathematics, Texas A&M University, College Station, TX.

Many local multiscale model reduction techniques have been developed previously. These approaches solve the underlying fine-scale problems on a coarse grid. Among these approaches are upscaling techniques [15, 38] and multiscale methods [6, 16, 21, 23, 24, 27, 7, 12, 8, 9]. In the latter, multiscale basis functions are locally constructed that capture local information. Many research papers [25, 37, 22] have been dedicated to optimizing limited number of multiscale basis functions to capture the solution accurately. In some recent works [18, 12, 9, 19, 20], the authors develop Generalized Multiscale Finite Element method (GMsFEM). GMsFEM is a flexible general framework that generalizes the Multiscale Finite Element Method (MsFEM) ([32]) by systematically enriching the coarse spaces. The main idea of this enrichment is to add extra basis functions that are needed to reduce the error substantially. This approach, as in many multiscale model reduction techniques, divides the computation into two stages: the offline and the online. In the offline stage, a small dimensional space is constructed that can be used in the online stage to construct multiscale basis functions. These multiscale basis functions can be re-used for any input parameter to solve the problem on a coarse grid. The main idea behind the construction of offline and online spaces is the selection of local spectral problems and the selection of the snapshot space.

In subsequent papers [10, 13], an adaptive GMsFEM is proposed. In these papers, we study an adaptive enrichment procedure and derive an a-posteriori error indicator which gives an estimate of the local error over coarse grid regions. The error indicators based on the $L^2$-norm of the local residual and on the weighted $H^{-1}$-norm of the local residual, where the weight is related to the coefficient of the elliptic equation are developed. We have shown that the use of weighted $H^{-1}$-norm residual gives a more robust error indicator which works well for cases with high contrast media. The error indicators contain the eigenvalue structure associated with GMsFEM. In particular, the smallest eigenvalue whose corresponding eigenvector is not included in the space enters into the error indicators.

Adaptivity is important for local multiscale methods as it identifies regions with large errors. However, after adding some initial basis functions, one needs to take into account some global information as the distant effects can be important. In this paper, we discuss the development of online basis functions that substantially accelerate the convergence of GMsFEM. The online basis functions are constructed based on a residual and motivated by the analysis.

We show, both theoretically and numerically, that one needs to have a sufficient number of initial basis functions in the offline space to guarantee an error decay independent of the contrast. We define such spaces as having online error reduction property (ONERP) and show that the eigenvalue that the corresponding eigenvector is not included in the offline space controls the error decay of the multiscale method. Larger is this eigenvalue, larger is the decrease in the error. Consequently, one needs to guarantee that eigenvectors associated with small (asymptotically small) eigenvalues are included in the initial coarse space. As we have discussed in [17, 26], many multiscale problems with high contrast can have very small eigenvalues and, thus, we need to include the eigenvectors associated with small eigenvalues in the initial coarse space.

Numerical results are presented to demonstrate that one needs to have a sufficient number of initial basis functions in the offline space before constructing online multiscale basis functions. Moreover, we study how different dimensional offline spaces can affect the error decay when online multiscale basis functions are added. We consider several examples where we vary the dimension of the offline space and add multiscale basis functions based on the residual. Our numerical results show that without sufficient number of offline basis functions, the error decay is not substantial. We study the proposed online basis construction in conjunction with adaptivity ([10, 13]), where online basis functions are added in some selected regions. Indeed, adaptivity is an important step to obtain an overall efficient local multiscale model reduction as it is essential to reduce the cost of online multiscale basis computations. Our numerical results show that the adaptive addition of online basis functions substantially improves GMsFEM. To reduce the computational cost associated with online multiscale basis computations, we propose computing the online basis functions in a reduced dimensional space consisting of several consequent offline basis functions. Our results show that one can still achieve a substantial error reduction this way. Because the online multiscale basis functions are not sparse in the offline space, approaches based on sparsity is not very helpful in our methods as our numerical results show.

In conclusion, the paper is organized in the following way. In Section 2, we present the underlying problem, the concepts of coarse and fine grids, and GMsFEM. In Section 3, we present some existing results for adaptive GMsFEM. In Section 4, we present our new proposed method for computing online multiscale basis functions. Numerical results are presented in Section 5. In Section 6, conclusions are drawn.

## 2 GMsFEM for high contrast flow

### 2.1 Overview

In this section, we will present a brief outline of the GMsFEM ([18, 12, 9, 19, 20]). Let $D$ be the computational domain. The high-contrast flow problem considered in this paper is

$$-\text{div}\big(\kappa(x)\nabla u\big) = f \quad \text{in} \quad D, \tag{1}$$

with the homogeneous Dirichlet boundary condition $u = 0$ on $\partial D$, where $f$ is a given source function. The difficulty in the numerical approximation of problem (1) arises from the complexity of the coefficient $\kappa(x)$, which can have multiple scales and very high contrast. In particular, discretizing (1) by traditional numerical schemes based on finite element or discontinuous Galerkin methods (e.g. [11, 35]) will result in very large and ill-conditioned linear systems, which require large computational times and memory. It is therefore desirable to develop efficient numerical schemes with a small number of degrees of freedom.

Next, we will introduce some notations. We use $\mathcal{T}^H$ to denote a usual conforming partition of the computational domain $D$. The set $\mathcal{T}^H$ is called the coarse grid and the elements of $\mathcal{T}^H$ are called coarse elements. Moreover, $H > 0$ is the coarse mesh size. In this paper, we consider rectangular coarse elements for the ease of discussions and illustrations. The methodology presented can be easily extended to coarse elements with more general geometries. Let $N$ be the number of nodes in the coarse grid $\mathcal{T}^H$, and let $\{x_i \, | \, 1 \leq i \leq N\}$ be the set of nodes in the coarse grid (or coarse nodes for short). For each coarse node $x_i$, we define a coarse neighborhood $\omega_i$ by

$$\omega_i = \bigcup \{K_j \in \mathcal{T}^H; \quad x_i \in \overline{K}_j\}. \tag{2}$$

Notice that $\omega_i$ is the union of all coarse elements $K_j \in \mathcal{T}^H$ having the coarse node $x_i$. An illustration of the above definitions is shown in Figure 1.

We let $\mathcal{T}^h$ be a partition of the computational domain $D$ obtained by refining the coarse grid $\mathcal{T}^H$. We call $\mathcal{T}^h$ the fine grid and $h > 0$ the fine mesh size. We remark that the restrictions of the fine grid in coarse neighborhoods will be used to discretize some local problems required for the generation of local basis functions. We will discuss this in the next section. Another use of the fine grid is for the computation of a fine-scale solution, which is used as a reference solution for comparison purposes. To fix the notations, we will use the standard conforming piecewise linear finite element method for the computation of the fine-scale solution. Specifically, we let $V$ be the conforming piecewise linear finite element space with respect to the fine grid $\mathcal{T}^h$. We will then obtain the fine-scale solution $u \in V$ by solving the following variational problem

$$a(u,v) = (f,v), \quad \text{for all } v \in V, \tag{3}$$

where $a(u,v) = \int_D \kappa(x) \nabla u \cdot \nabla v \, dx$, and $(f,v) = \int_D fv \, dx$. Note that we equip the space $V$ with the energy norm $\|v\|_V^2 = a(v,v)$. We assume that the fine mesh size $h$ is small enough so that the fine-scale solution $u$ is close enough to the exact solution. The purpose of this paper is to find a multiscale solution $u_{\text{ms}}$ that is a good approximation of the fine-scale solution $u$. The multiscale solution $u_{\text{ms}}$ is obtained by the GMsFEM.

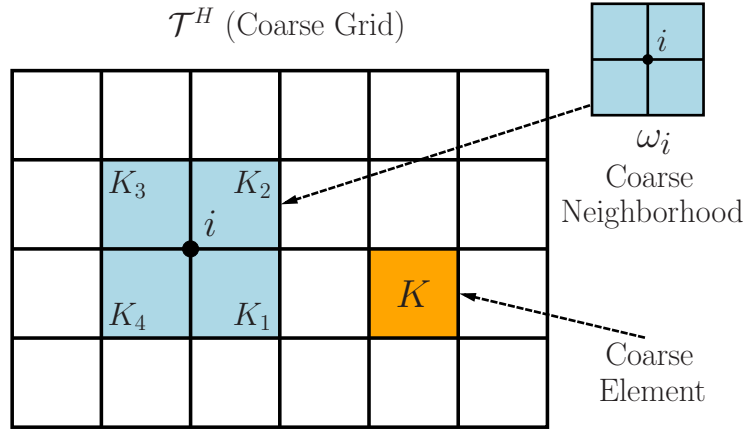

Figure 1: Illustration of a coarse neighborhood and a coarse element.

Now we present the general idea of GMsFEM ([18, 12, 9, 19, 20]). We will consider the continuous Galerkin (CG) formulation, which has a similar form as the fine-scale problem (3). The basis functions are nodal based and have supports on coarse neighborhoods. Specifically, for each coarse node $x_i$, we will construct a set of basis functions $\{\psi_k^{\omega_i} \mid k = 1, 2, \cdots, l_i\}$ such that each $\psi_k^{\omega_i}$ is supported on the coarse neighborhood $\omega_i$, where $l_i$ is the number of basis functions with support in $\omega_i$. In addition, the basis functions satisfy a partition of unity property, namely, there are coefficients $\alpha_k^i$ such that $\sum_{i=1}^{N} \sum_{k=1}^{l_i} \alpha_k^i \psi_k^{\omega_i} = 1$. We remark that contrary to standard multiscale finite element method, one can use multiple basis functions for each coarse node and use different numbers of basis functions for different coarse nodes in our GMsFEM. Once the basis functions are constructed, we can define the approximation space $V_{\text{ms}}$ by the linear span of all basis functions. The GMsFEM solution $u_{\text{ms}} \in V_{\text{ms}}$ can then be obtained by solving the following

$$a(u_{\text{ms}}, v) = (f,v), \quad \text{for all } v \in V_{\text{ms}}. \tag{4}$$

We remark that one can also use the discontinuous Galerkin formulation (see e.g., [8, 9, 19]) instead of the CG formulation.

From the above, one sees that the key ingredient of the GMsFEM is the construction of local basis functions. Based on the works [18, 12, 9, 19, 20], we will use the so called offline basis functions, which can be computed in the offline stage. Moreover, we will construct online basis functions that are problem dependent and are computed locally and adaptively based on some local residuals. Our results show that the combination of both offline and online basis functions will give a rapid convergence of the multiscale solution $u_{\text{ms}}$ to the fine-scale solution $u$.

## 2.2 Construction of offline basis functions

In this section, we will present the construction of the offline basis functions (see e.g. [18, 12, 9, 19, 20]). Let $\omega$ be a given coarse neighborhood. Notice that we omit the coarse node index to simplify the notations. The construction begins with a snapshot space $V_{\text{snap}}^{\omega}$. The snapshot space $V_{\text{snap}}^{\omega}$ is a set of functions defined on $\omega$ and contains all or most necessary components of the fine-scale solution restricted to $\omega$. A spectral problem is then solved in the snapshot space to extract the dominant modes in the snapshot space. These dominant modes are the offline basis functions and the resulting reduced space is called the offline space. There are two choices of $V_{\text{snap}}^{\omega}$ that are commonly used. The first choice is the restriction of the conforming space $V$ in $\omega$, and the resulting basis functions are called **spectral basis functions**. The second choice is the set of all $\kappa$-harmonic extensions, and the resulting basis functions are called **harmonic basis functions**.

Next, we recall the definition of the snapshot space $V_{\text{snap}}^{\omega}$ based on $\kappa$-harmonic extensions. Let $J_h(\omega_i)$ be the set of all nodes of the fine mesh $\mathcal{T}^h$ lying on $\partial \omega_i$. For each fine-grid node $x_j \in J_h(\omega_i)$, we define a discrete delta function $\delta_j^h(x)$ defined in $J_h(\omega_i)$ by

$$\delta_j^h(x_k) = \begin{cases} 1, & k = j \\ 0, & k \neq j \end{cases}, \quad x_k \in J_h(\omega_i).$$

The $j$-th snapshot function $\psi_j^{\omega,\text{snap}} \in V_{\text{snap}}^{\omega}$ for the coarse neighborhood $\omega$ is defined as the solution of

$$\begin{aligned} -\text{div}(\kappa(x) \nabla \psi_j^{\omega,\text{snap}}) &= 0, \quad \text{in } \omega, \\ \psi_j^{\omega,\text{snap}} &= \delta_j^h, \quad \text{on } \partial \omega. \end{aligned} \quad (5)$$

Clearly, the dimension of $V_{\text{snap}}^{\omega}$ is equal to the number of elements in $J_h(\omega_i)$, the set of fine-grid nodes lying on $\partial \omega$. We note that one can use randomized snapshots in conjunction with oversampling to reduce the computational cost associated with the snapshot calculations. We refer to [5] for details.

To obtain the offline basis functions, we need to perform a space reduction by a spectral problem. The analysis in [21] motivates the following construction. The spectral problem that is needed for the purpose of space reduction is: find $(\psi, \lambda) \in V_{\text{snap}}^{\omega} \times \mathbb{R}$ such that

$$\int_{\omega} \kappa(x) \nabla \psi \cdot \nabla \phi \, dx = \lambda \int_{\omega} \widetilde{\kappa}(x) \psi \phi \, dx, \quad \forall \phi \in V_{\text{snap}}^{\omega} \quad (6)$$

where the weighted function $\widetilde{\kappa}(x)$ is defined by (see [21])

$$\widetilde{\kappa} = \kappa \sum_{i=1}^{N} H^2 |\nabla \chi_i|^2,$$

and $\chi_i$ is the standard multiscale basis function for the coarse node $x_i$ (that is, with linear boundary conditions for cell problems). More precisely,

$$-\text{div}\,(\kappa(x)\nabla\chi_i) = 0 \quad K \in \omega_i \tag{7}$$
$$\chi_i = g_i \quad \text{on } \partial K,$$

for all $K \in \omega_i$, where $g_i$ is a continuous function on $\partial K$ and is linear on each edge of $\partial K$. We arrange the eigenvalues $\lambda_k^\omega$, $k = 1, 2, \cdots$, from (6) in ascending order. We then select the first $l_i$ eigenfunctions from (6), and denote them by $\Psi_1^{\omega,\text{off}}, \cdots, \Psi_{l_i}^{\omega,\text{off}}$. Using these eigenfunctions, we can define

$$\phi_k^{\omega,\text{off}} = \sum_{j=1}^{l_i} (\Psi_k^{\omega,\text{off}})_j \psi_j^{\omega,\text{snap}}, \qquad k = 1, 2, \cdots, l_i,$$

where $(\Psi_k^{\omega,\text{off}})_j$ denotes the $j$-th component of $\Psi_k^{\omega,\text{off}}$. Finally, the offline basis functions for the coarse neighborhood $\omega$ is defined by $\psi_k^{\omega,\text{off}} = \chi \phi_k^{\omega,\text{off}}$, where $\chi$ is the standard multiscale basis function for a generic coarse neighborhood $\omega$. We also define the local offline space $V_{\text{off}}^\omega$ as the linear span of all $\psi_k^{\omega,\text{off}}$, $k = 1, 2, \cdots, l_i$.

We remark that one can take $V_{\text{ms}}$ in (4) as $V_{\text{off}} := \text{span}\{\psi_k^{\omega_i,\text{off}} \mid 1 \leq i \leq N, 1 \leq k \leq l_i\}$. The convergence of the resulting method is analyzed in [21].

## 3 Offline Adaptive GMsFEM

The use of $V_{\text{ms}} = V_{\text{off}}$ in (4) is a promising option in a variety of scenarios. However, the number of basis functions $l_i$ used for the coarse neighborhood $\omega_i$ has to be pre-defined in [21]. Recently, an adaptive enrichment algorithm is proposed and analyzed in [10], allowing $l_i$ to be chosen adaptively. The algorithm allows one to use more basis functions in regions with more complexity without using a priori information. To be more specific, we consider a coarse neighborhood $\omega_i$ and assume that $l_i$ basis functions are currently used. One can then compute a residual based on the current solution. If the residual is large according to a certain criteria, we will add one (or more) basis function(s) by using the next eigenfunction(s) $\Psi_{l_i+1}^{\omega_i,\text{off}}, \Psi_{l_i+2}^{\omega_i,\text{off}}, \cdots$. This iterative process is stopped when some error tolerance is reached. Furthermore, the convergence of this algorithm is proved in [10], with a convergence rate independent of contrasts in $\kappa(x)$. We remark that we call this algorithm the offline adaptive GMsFEM since only offline basis functions are used. In the next section, we will discuss the construction and the use of online basis functions.

Now, we will briefly state the offline adaptive GMsFEM. Let $u_{\text{ms}} \in V_{\text{off}}$ be the solution obtained in (4). Consider a given coarse neighborhood $\omega_i$. We define a space $V_i = H_0^1(\omega_i) \cap V$ which is equipped with the norm $\|v\|_{V_i}^2 = \int_{\omega_i} \kappa(x)|\nabla v|^2\,dx$. We also define the following linear functional on $V_i$ by

$$R_i(v) = \int_{\omega_i} fv - \int_{\omega_i} a\nabla u_{\text{ms}} \cdot \nabla v. \tag{8}$$

This is called the $H^{-1}$-residual on $\omega_i$. The functional norm of $R_i$, denoted by $\|R_i\|_{V_i^*}$, gives a measure of the size of the residual. The first important result in [10] states that these residuals give a computable indicator of the error $u - u_{\text{ms}}$ in the energy norm. In particular, we have

$$\|u - u_{\text{ms}}\|_V^2 \leq C_{\text{err}} \sum_{i=1}^{N} \|R_i\|_{V_i^*}^2 (\lambda_{l_i+1}^{\omega_i})^{-1}, \tag{9}$$

where $C_{\text{err}}$ is a uniform constant, and $\lambda_{l_i+1}^{\omega_i}$ denotes the $(l_i+1)$-th eigenvalue for the problem (6) in the coarse neighborhood $\omega_i$, and corresponds to the first eigenvector that is not included in the construction of $V_{\text{off}}^{\omega}$.

The adaptive enrichment algorithm [10] is stated as follows. We use the index $m \geq 1$ to represent the enrichment level. At the enrichment level $m$, we use $V_{\text{off}}^m$ to denote the corresponding GMsFEM space and $u_{\text{ms}}^m$ the corresponding solution obtained in (4), with $V_{\text{ms}} = V_{\text{off}}^m$. Furthermore, we use $l_i^m$ to denote the number of basis functions used in the coarse neighborhood $\omega_i$. We will present the strategy for getting the space $V_{\text{off}}^{m+1}$ from $V_{\text{off}}^m$. Let $0 < \theta < 1$ be a given number independent of $m$. First of all, we compute the local residuals for every coarse neighborhood $\omega_i$:

$$\eta_i^2 = \|R_i\|_{V_i^*}^2 (\lambda_{l_i^m+1}^{\omega_i})^{-1}, \qquad i = 1, 2, \cdots, N,$$

where $R_i(v)$ is defined using $u_{\text{ms}}^m$, namely,

$$R_i(v) = \int_{\omega_i} fv - \int_{\omega_i} a\nabla u_{\text{ms}}^m \cdot \nabla v, \quad \forall v \in V_i.$$

Next, we will add basis functions for the coarse neighborhoods with large residuals. To do so, we re-enumerate the coarse neighborhoods so that the above local residuals $\eta_i^2$ are arranged in decreasing order $\eta_1^2 \geq \eta_2^2 \geq \cdots \geq \eta_N^2$. We then select the smallest integer $k$ such that

$$\theta \sum_{i=1}^{N} \eta_i^2 \leq \sum_{i=1}^{k} \eta_i^2. \tag{10}$$

For those coarse neighborhoods $\omega_1, \cdots, \omega_k$ (in the new enumeration) chosen in the above procedure, we will add basis functions by using the next eigenfunctions $\Psi_{l_i+1}^{\omega_i,\text{off}}, \Psi_{l_i+2}^{\omega_i,\text{off}}, \cdots$. The resulting space is called $V_{\text{off}}^{m+1}$. We remark that the choice of $k$ defined in (10) is called the Dorlfer's bulk marking strategy [14]. For more details about this enrichment algorithm, see [10].

Both numerical and theoretical results in [10] show that the enrichment algorithm gives a rapid convergence when the eigenvalues in (6) have a fast growth. While there is a large class of problems having fast eigenvalue growth, there are still cases for which the eigenfunctions corresponding to large eigenvalues have little contribution to the fine-scale solution. One reason is that the local basis functions do not contain any global information, and thus they cannot be used to efficiently capture these global behaviors. We will therefore present in the next section that some online basis functions are necessary to obtain a coarse representation of the fine-scale solution and give a rapid convergence of the corresponding adaptive enrichment algorithm.

## 4 Residual based online adaptive GMsFEM

As we mentioned in the previous sections, some online basis functions are necessary to obtain a coarse representation of the fine-scale solution and give a rapid convergence of the corresponding adaptive enrichment algorithm. In this section, we will give the precise meaning of online basis functions and the corresponding adaptive enrichment algorithm. We will first derive a framework for the constructions of online multiscale basis functions. Based on our derivations, we will argue that one also needs offline basis functions to satisfy some properties in order to guarantee that adding online basis functions will decrease the error.

We will use similar notations as in the previous section. We use the index $m \geq 1$ to represent the enrichment level. At the enrichment level $m$, we use $V_{\text{ms}}^m$ to denote the corresponding GMsFEM space and $u_{\text{ms}}^m$ the corresponding solution obtained in (4). The sequence of functions $\{u_{\text{ms}}^m\}_{m\geq 1}$ will converge to the fine-scale solution. We emphasize that the space $V_{\text{ms}}^m$ can contain both offline and online basis functions. We will construct a strategy for getting the space $V_{\text{ms}}^{m+1}$ from $V_{\text{ms}}^m$.

Next we present a framework for the construction of online basis functions. By online basis functions, we mean basis functions that are computed during the iterative process, contrary to offline basis functions that are computed before the iterative process. The online basis functions are computed based on some local residuals for the current multiscale solution, that is, the function $u_{\text{ms}}^m$. Thus, we see that some offline basis functions are necessary for the computations of online basis functions. We will also see how many of these offline basis functions are needed in order to obtain a rapidly converging sequence of solutions.

Consider a given coarse neighbourhood $\omega_i$. Suppose that we need to add a basis function $\phi \in V_i$ on the $i$-th coarse neighbourhood $\omega_i$. Let $V_{\text{ms}}^{m+1} = V_{\text{ms}}^m + \text{span}\{\phi\}$ be the new approximation space, and $u_{\text{ms}}^{m+1} \in V_{\text{ms}}^{m+1}$ be the corresponding GMsFEM solution. It is easy to see from (4) that $u_{\text{ms}}^{m+1}$ satisfies

$$\|u - u_{\text{ms}}^{m+1}\|_V^2 = \inf_{v \in V_{\text{ms}}^{m+1}} \|u - v\|_V^2.$$

Taking $v = u_{\text{ms}}^m + \alpha \phi$, where $\alpha$ is a scalar to be determined, we have

$$\|u - u_{\text{ms}}^{m+1}\|_V^2 \leq \|u - u_{\text{ms}}^m - \alpha \phi\|_V^2 = \|u - u_{\text{ms}}^m\|_V^2 - 2\alpha a(u - u_{\text{ms}}^m, \phi) + \alpha^2 a(\phi, \phi).$$

The last two terms in the above inequality measure the amount of reduction in error when the new basis function $\phi$ is added to the space $V_{\text{ms}}^m$. To determine $\phi$, we first assume that the basis function $\phi$ is normalized so that $a(\phi, \phi) = 1$. In order to maximize the reduction in error, we will find $\alpha$ in order to maximize the quantity $2\alpha a(u - u_{\text{ms}}^m, \phi) - \alpha^2$. Clearly, one needs to take $\alpha = a(u - u_{\text{ms}}^m, \phi)$. Using this choice of $\alpha$, we have

$$\|u - u_{\text{ms}}^{m+1}\|_V^2 \leq \|u - u_{\text{ms}}^m\|_V^2 - |a(u - u_{\text{ms}}^m, \phi)|^2.$$

Since $\phi \in V_i \subset V$, by using (3), we have

$$\|u - u_{\text{ms}}^{m+1}\|_V^2 \leq \|u - u_{\text{ms}}^m\|_V^2 - |(f, \phi) - a(u_{\text{ms}}^m, \phi)|^2.$$

We will then find $\phi \in V_i$ to maximize the local residual $|(f, \phi) - a(u_{\text{ms}}^m, \phi)|^2$. Clearly, the maximum of the quantity $|(f, \phi) - a(u_{\text{ms}}^m, \phi)|$ equals to the functional norm of the residual $R_i$. Moreover, the required $\phi \in V_i$ is the solution of

$$a(\phi, v) = (f, v) - a(u_{\text{ms}}^m, v), \quad \forall v \in V_i \tag{11}$$

and $\|\phi\|_{V_i} = \|R_i\|_{V_i^*}$. Hence, the new online basis function $\phi \in V_i$ can be obtained by solving (11). In addition, the residual norm $\|R_i\|_{V_i^*}$ provides a measure on the amount of reduction in energy error. We remark that we call this algorithm the online adaptive GMsFEM since only online basis functions are used.

Now, we study the convergence of the above online adaptive procedure. To simplify notations, we write $r_i = \|R_i\|_{V_i^*}$. From the above constructions, we have

$$\|u - u_{\text{ms}}^{m+1}\|_V^2 \leq \|u - u_{\text{ms}}^m\|_V^2 - r_i^2. \tag{12}$$

We assume that each of the spaces $V_{\text{ms}}^m$, $m \geq 1$, contains $n_j$ offline basis functions for the coarse neighborhood $\omega_j$. Then, similar to (9), we have

$$\|u - u_{\text{ms}}^m\|_V^2 \leq C_{\text{err}} \sum_{j=1}^N r_j^2 (\lambda_{n_j+1}^{\omega_j})^{-1}. \tag{13}$$

Combining (12) and (13), we obtain

$$\|u - u_{\text{ms}}^{m+1}\|_V^2 \leq \left(1 - \frac{\lambda_{n_i+1}^{\omega_i}}{C_{\text{err}}} \frac{r_i^2(\lambda_{n_i+1}^{\omega_i})^{-1}}{\sum_{j=1}^N r_j^2(\lambda_{n_j+1}^{\omega_j})^{-1}}\right) \|u - u_{\text{ms}}^m\|_V^2.$$

The above inequality gives the convergence of the online adaptive GMsFEM with a precise convergence rate for the case when one online basis function is added per iteration. To enhance the convergence and efficiency of the online adaptive GMsFEM, we consider enrichment on non-overlapping coarse neighborhoods. Let $I \subset \{1, 2, \cdots, N\}$ be the index set of some non-overlapping coarse neighborhoods. For each $i \in I$, we can obtain a basis function $\phi_i \in V_i$ using (11). We define $V_{\text{ms}}^{m+1} = V_{\text{ms}}^m + \text{span}\{\phi_i, i \in I\}$. Following the same argument as above and using the fact that the coarse neighborhoods $\omega_i$, $i \in I$, are non-overlapping, we obtain

$$\|u - u_{\text{ms}}^{m+1}\|_V^2 \leq \|u - u_{\text{ms}}^m\|_V^2 - \sum_{i \in I} r_i^2. \tag{14}$$

Consequently, we have

$$\|u - u_{\text{ms}}^{m+1}\|_V^2 \leq \left(1 - \frac{\Lambda_{\min}^{(I)}}{C_{\text{err}}} \frac{\sum_{i \in I} r_i^2(\lambda_{n_i+1}^{\omega_i})^{-1}}{\sum_{j=1}^N r_j^2(\lambda_{n_j+1}^{\omega_j})^{-1}}\right) \|u - u_{\text{ms}}^m\|_V^2 \tag{15}$$

where

$$\Lambda_{\min}^{(I)} = \min_{i \in I} \lambda_{n_i+1}^{\omega_i}.$$

Inequality (15) shows that we are able to obtain a better convergence of our online adaptive GMsFEM by adding more online basis functions per iteration. The convergence rate depends on the factors $C_{\text{err}}$ and $\Lambda_{\min}^{(I)}$. We will therefore need to take enough offline basis functions so that $\Lambda_{\min}^{(I)}$ is large enough and

$$\frac{\Lambda_{\min}^{(I)}}{C_{\text{err}}} \frac{\sum_{i \in I} r_i^2(\lambda_{n_i+1}^{\omega_i})^{-1}}{\sum_{j=1}^N r_j^2(\lambda_{n_j+1}^{\omega_j})^{-1}} \geq \theta_0$$

for some $0 < \theta_0 < 1$ which is independent of the contrast in $\kappa(x)$. Hence, we obtain the following convergence for the online adaptive GMsFEM:

$$\|u - u_{\text{ms}}^{m+1}\|_V^2 \leq (1 - \theta_0) \|u - u_{\text{ms}}^m\|_V^2.$$

We note that $\Lambda_{\min}^{(I)}$ can be very small when there are channels in the domain. This is extensively discussed in [17]. For this reason, we introduce a definition.

**Definition 4.1.** *We say $V_{\text{off}}$ satisfies Online Error Reduction Property (ONERP) if*

$$\frac{\Lambda_{\min}^{(I)}}{C_{\text{err}}} \frac{\sum_{i \in I} r_i^2(\lambda_{n_i+1}^{\omega_i})^{-1}}{\sum_{i=1}^N r_i^2(\lambda_{n_i+1}^{\omega_i})^{-1}} \geq \theta_0,$$

*for some $\theta_0 > \delta > 0$, where $\delta$ is independent of physical parameters such as contrast.*

We remark that if $V_{\text{off}}$ is ONERP, then the error will decrease independent of physical parameters such as the contrast and scales. We will show in our numerical results that if we do not choose $V_{\text{off}}$ with ONERP, the online basis functions will not decrease the error. One of easiest way to determine $V_{\text{off}}$ being ONERP is to guarantee that $\Lambda_{\min}^{(I)}$ is sufficiently large. In general, one can use the sizes of $\Lambda_{\min}^{(I)}$ and $\sum_{i \in I} r_i^2 (\lambda_{n_i+1}^{\omega_i})^{-1}$ to determine the switching between offline and online.

## 5 Numerical result

In this section, we will present numerical examples to demonstrate the performance of the proposed method. The online adaptive GMsFEM is implemented as follows. We will first choose a fixed number of offline basis functions for every coarse neighborhood, and denote the resulting offline space as $V_{\text{off}}$. We set $V_{\text{ms}}^1 = V_{\text{off}}$. In addition, the basis functions in $V_{\text{ms}}^1$ are called *initial basis*. We will enumerate the coarse neighborhoods by a two-index notation. More precisely, the coarse neighborhoods are denoted by $\omega_{i,j}$, where $i = 1, 2, \cdots, N_x$ and $j = 1, 2, \cdots, N_y$ and $N_x$ and $N_y$ are the number of coarse nodes in the $x$ and $y$ directions respectively. We let $I_{x,\text{odd}}$ and $I_{x,\text{even}}$ be the set of odd and even indices from $\{1, 2, \cdots, N_x\}$. We use similar definitions for $I_{y,\text{odd}}$ and $I_{y,\text{even}}$. Each iteration of our online adaptive GMsFEM contains 4 sub-iterations. In particular, these 4 sub-iterations are defined by adding online basis functions in the non-overlapping coarse neighborhoods $\omega_{i,j}$ with $(i,j) \in I_{x,\text{odd}} \times I_{y,\text{odd}}$, $(i,j) \in I_{x,\text{odd}} \times I_{y,\text{even}}$, $(i,j) \in I_{x,\text{even}} \times I_{y,\text{odd}}$ and $(i,j) \in I_{x,\text{even}} \times I_{y,\text{even}}$ respectively.

The domain $D$ is taken as the unit square $[0,1]^2$ and is divided into $16 \times 16$ coarse blocks consisting of uniform squares. Each coarse block is then divided into $16 \times 16$ fine blocks consisting of uniform squares. That is, the whole domain is partitioned by $256 \times 256$ fine grid blocks. The medium parameter $\kappa$ and the source function $f$ are shown in Figure 2. We will use the following

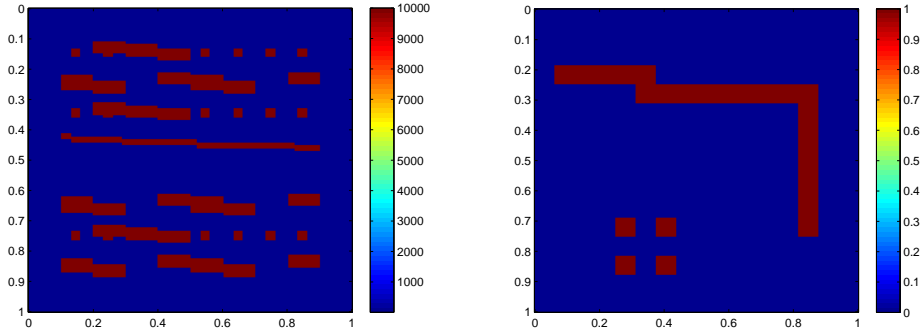

Figure 2: Left: Permeability field $\kappa$. Right: Source function $f$.

error quantities to compare the accuracy of our algorithm

$$e_2 = \frac{\|u - u_{\text{ms}}\|_{L^2(D)}}{\|u\|_{L^2(D)}}, \quad e_a = \frac{\|u - u_{\text{ms}}\|_V}{\|u\|_V}.$$

## 5.1 Comparison of using different number of initial basis

In Table 1, we present the convergence history of our algorithm for using one initial basis per coarse neighborhood. Notice that we have shown the number of basis functions used for each coarse neighborhood and the total degrees of freedom (DOF), which are the numbers in parentheses. We use the multiscale basis functions as the initial partition of unity (see (7)). We consider two different contrasts. On the right table, we increase the contrast by 100 times. More precisely, the conductivity of inclusions and channels in Figure 2 (left figure) is multiplied by 100. In this case, first few eigenvalues that are in the regions with channels become 100 times smaller ([17]). This decrease in the eigenvalues will make the error decay slower. This can be observed by comparing the left and the right tables of Table 1, where we can see that the errors in case of the higher contrast decrease much slower. This is also observed when we use 2 initial basis functions (see Table 2). When using two initial basis functions, there are contrast-dependent small eigenvalues (these eigenvalues decrease as we increase the contrast), and thus, by increasing the contrast, the decay becomes slower. This can be observed in Table 2. However, if we choose 3 initial basis functions, then $\Lambda_{\min}$ is independent of the contrast and, thus, for larger contrasts, we observe a similar error behavior. We observed a similar convergence when using 4 initial basis functions, see Figure 3, where we plot the relative energy error against the dimension of $V_{\text{ms}}$ for various choices of the initial basis and for two types of contrasts, $1e4$ (left) and $1e6$ (right).

| number of basis (DOF) | $e_a$ | $e_2$ | number of basis (DOF) | $e_a$ | $e_2$ |
|---|---|---|---|---|---|
| 1(225) | 60.71% | 33.87% | 1(225) | 60.90% | 34.15% |
| 2(450) | 33.01% | 13.38% | 2(450) | 35.90% | 15.87% |
| 3(675) | 14.38% | 3.25% | 3(675) | 35.00% | 15.29% |
| 4(900) | 4.28% | 1.02% | 4(900) | 25.77% | 8.77% |
| 5(1125) | 1.33% | 0.24% | 5(1125) | 14.17% | 4.39% |
| 6(1350) | 0.065% | 0.0028% | 6(1350) | 7.79% | 2.78% |
| 7(1575) | 0.00083% | 2.96e-05% | 7(1575) | 6.83% | 2.06% |
| 8(1800) | 1.59e-05% | 4.87e-07% | 8(1800) | 4.15% | 1.20% |
| 9(2025) | 2.35e-07% | 2.10e-08% | 9(2025) | 2.60% | 0.64% |

Table 1: One initial basis. Left: Lower contrast($1e4$). Right: Higher contrast($1e6$).

| num of basis(DOF) | $e_a$ | $e_2$ | num of basis (DOF) | $e_a$ | $e_2$ |
|---|---|---|---|---|---|
| 2 (450) | 26.60% | 6.92% | 2 (450) | 27.17% | 7.53% |
| 3 (675) | 1.46% | 0.060% | 3 (675) | 4.99% | 0.79% |
| 4 (900) | 0.017% | 0.000079% | 4 (900) | 0.20% | 0.0073% |
| 5 (1125) | 0.000021% | 1.06e-05% | 5 (1125) | 0.0017% | 8.16e-05 |
| 6 (1350) | 3.56e-06% | 1.65e-07% | 6 (1350) | 2.71e-05% | 1.09e-06% |

Table 2: Two initial basis. Left: Lower contrast($1e4$). Right: Higher contrast($1e6$).

Next, we will present an example with a different medium parameter $\kappa$ shown in Figure 4 to show the importance of ONERP. The source function $f$ is taken as the constant 1. The domain $D$ is divided into $8 \times 8$ coarse blocks consisting of uniform squares. Each coarse block is then divided into $32 \times 32$ fine blocks also consisting of uniform squares. The convergence history for the use of one initial basis and the corresponding total number of degrees of freedom (DOF) are shown in Table 4. In this case, $\Lambda_{\min} = 0.0033$, which is considered to be very small, and we observe very slow convergence of the online adaptive procedure. In Table 5, we present the convergence history for the use of two to five initial basis, where we only show the results for the last 4 iterations. We see that the values of $\Lambda_{\min}$ increase as we increase the number of initial basis. We also observe that the convergence rate increase when we raise the number of initial basis from 2 to 4. For the use of 5 initial basis, we again see rapid convergence and a faster convergence compared when using 4 initial basis functions. In particular, we observe (based on 3 iterations following the initial one) that the error decays at 130-fold when 5 initial basis functions are selected, while the error decay is about 90-fold when 4 initial basis functions are selected. A comparison of error decay for the use of 1 to 5 initial basis functions is shown in Figure 5. We have also tested harmonic basis functions and the results are similar, i.e., the convergence rate is very slow unless sufficient number of offline basis functions is selected.

| num of basis (DOF) | $e_a$ | $e_2$ | num of basis (DOF) | $e_a$ | $e_2$ |
|---|---|---|---|---|---|
| 3 (675) | 16.95% | 2.53% | 3 (675) | 16.96% | 2.54% |
| 4 (900) | 0.54% | 0.023% | 4 (900) | 0.54% | 0.023% |
| 5 (1125) | 0.011% | 0.00040% | 5 (1125) | 0.011% | 0.00041% |
| 6 (1350) | 9.07e-05% | 3.79e-06% | 6 (1350) | 9.07e-05% | 3.79e-06% |
| 7 (1575) | 1.38e-06% | 6.05e-08% | 7 (1575) | 1.58e-06% | 5.49e-07% |

Table 3: Three initial basis. Left: Lower contrast(1$e$4). Right: Higher contrast(1$e$6).

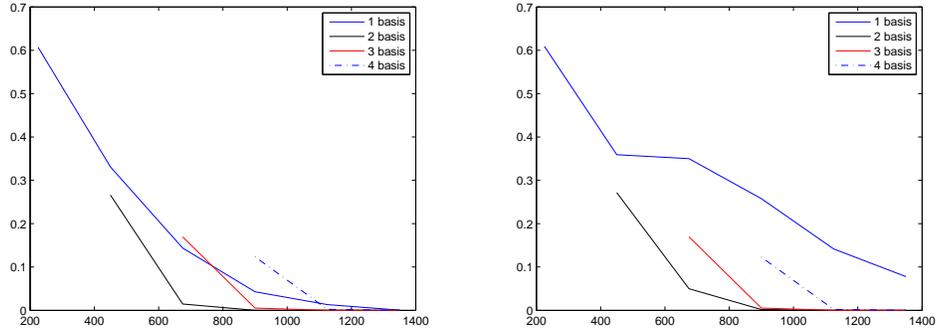

Figure 3: Error comparison. Along $x$-axis: Dimensions of $V_{\text{ms}}$. Along $y$-axis: Relative energy errors. Left: 1$e$4. Right: 1$e$6.

In conclusion, we observe

- If $V_{\text{ms}}$ does not satisfy ONERP, then the error decay is slower as the contrast becomes larger.

- If $V_{\text{ms}}$ does not satisfy ONERP, in some cases, we have observed the error does not decrease as we add online basis functions (see Table 4, 5).

- If $V_{\text{ms}}$ satisfies ONERP, then we observe a fast convergence, which is independent of contrast.

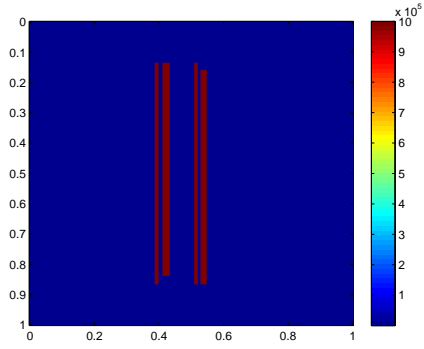

Figure 4: Permeability field $\kappa$.

| DOF | $e_a$ | $e_2$ |
|---|---|---|
| 81 | 17.24% | 4.35% |
| 162 | 2.80% | 1.02% |
| 243 | 2.65% | 0.88% |
| 323 | 2.64% | 0.87% |
| 401 | 1.09% | 0.081% |
| 478 | 0.74% | 0.094% |
| 555 | 0.73% | 0.090% |
| 632 | 0.48% | 0.039% |
| 709 | 0.37% | 0.026% |

Table 4: One initial basis ($\Lambda_{\min} = 0.0033$).

| DOF | $e_a$ | $e_2$ |
|---|---|---|
| 162 | 13.29% | 2.90% |
| 243 | 2.00% | 0.32% |
| 324 | 1.79% | 0.23% |
| 405 | 1.60% | 0.17% |
| 486 | 0.33% | 0.025% |
| 567 | 0.30% | 0.022% |
| 648 | 0.012% | 0.00057% |
| 725 | 0.00012% | 4.99e-06% |
| 765 | 3.62e-06% | 2.45e-06% |

| DOF | $e_a$ | $e_2$ |
|---|---|---|
| 243 | 10.26% | 1.51% |
| 324 | 1.78% | 0.23% |
| 405 | 1.75% | 0.23% |
| 486 | 0.25% | 0.0088% |
| 567 | 0.0016% | 0.000089% |
| 638 | 8.34e-06% | 2.50e-06% |
| 644 | 3.69e-06% | 2.49e-06% |

| DOF | $e_a$ | $e_2$ |
|---|---|---|
| 324 | 7.95% | 1.06% |
| 405 | 0.074% | 0.0035% |
| 486 | 0.0010% | 3.24e-05% |
| 563 | 1.10e-05% | 2.55e-06% |
| 568 | 3.71e-06% | 2.52e-06% |

| DOF | $e_a$ | $e_2$ |
|---|---|---|
| 405 | 7.24% | 0.92% |
| 486 | 0.0684% | 0.0028% |
| 567 | 0.00049% | 1.51e-05% |
| 635 | 3.80e-06% | 2.49e-06% |

Table 5: Top-left: Two initial basis ($\Lambda_{\min} = 0.026$). Top-right: Three initial basis ($\Lambda_{\min} = 0.18$), Bottom-left: Four initial basis ($\Lambda_{\min} = 199.12$). Bottom-right: Five initial basis ($\Lambda_{\min} = 319.32$).

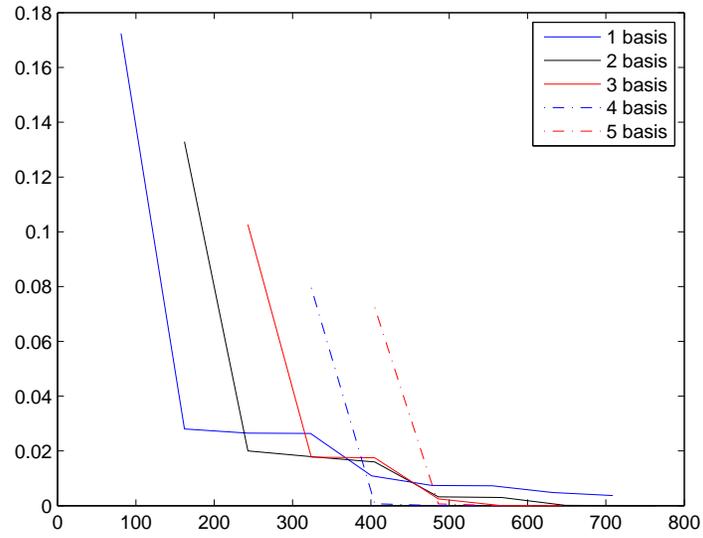

Figure 5: Error comparison for different number of initial basis functions. Along $x$-axis: Dimensions of $V_{\text{ms}}$. Along $y$-axis: Relative energy errors.

## 5.2 Adaptive online enrichment

In this section, the online enrichment is performed only for regions with the residual that is larger than a certain threshold. In the first case, the online enrichment is performed for the coarse regions with a residual error bigger than a certain threshold which will be taken $10^{-3}$, $10^{-4}$, and $10^{-5}$. In the second case, the online enrichment is performed for coarse regions that have cumulative residual that is $\theta$ fraction of the total residual. One of our objectives is to show that one can drive the error down to a number below a threshold, adaptively.

In our numerical results, we will consider three tolerances ($tol$) $10^{-3}$, $10^{-4}$ and $10^{-5}$. We will enrich coarse regions, if the $H^{-1}$-norm of the residual is bigger than the tolerance. In Table 6, we show the errors when using 1 initial basis function for tolerances $10^{-3}$, $10^{-4}$ and $10^{-5}$. We first observe a very slow reduction in errors similar to the results presented in the previous section. Another observation is that the energy error of the multiscale solution is in the same order of the tolerance, and the error cannot be further reduced if we perform more iterations. This allows us to compute a multiscale solution with a prescribed error level by choosing a suitable tolerance in the adaptive algorithm. In Table 7 and Table 8, we show the errors for the last three iterations when using 2 and 3 initial basis functions respectively for tolerances $10^{-3}$, $10^{-4}$ and $10^{-5}$. We observe that the convergences are much faster. In addition, the energy errors are again have the same magnitude as the tolerances. From these results, we obtain the following conclusions.

- Using smaller tolerances, we can reduce the final error below desired threshold errors.

- We have observed that the number of initial basis functions are important to achieve better results. For example, we observe a slow decay of the error when 1 initial basis function is selected. Moreover, if the contrast is higher, the decay becomes slower.

| DOF | $e_a$ | $e_2$ | DOF | $e_a$ | $e_2$ | DOF | $e_a$ | $e_2$ |
|---|---|---|---|---|---|---|---|---|
| 225 | 60.71% | 33.87% | 225 | 60.71% | 33.87% | 225 | 60.701% | 33.87% |
| 447 | 33.10% | 13.39% | 449 | 33.10% | 13.38% | 450 | 33.10% | 13.38% |
| 652 | 14.43% | 3.28% | 674 | 14.38% | 3.25% | 675 | 14.38% | 3.25% |
| 776 | 4.37% | 1.06% | 883 | 4.28% | 1.02% | 899 | 4.28% | 1.02% |
| 824 | 1.83% | 0.37% | 1031 | 1.33% | 0.24% | 1114 | 1.33% | 0.24% |
| 847 | 1.10% | 0.25% | 1125 | 0.082% | 0.0036% | 1275 | 0.065% | 0.0028% |
| 863 | 0.50% | 0.029% | 1136 | 0.052% | 0.0023% | 1338 | 0.0048% | 0.00017% |

Table 6: One initial basis. Left: $tol = 10^{-3}$. Middle: $tol = 10^{-4}$. Right: $tol = 10^{-5}$.

| DOF | $e_a$ | $e_2$ | DOF | $e_a$ | $e_2$ | DOF | $e_a$ | $e_2$ |
|---|---|---|---|---|---|---|---|---|
| 450 | 26.60% | 6.92% | 450 | 26.60% | 6.92% | 675 | 1.46% | 0.060% |
| 649 | 1.49% | 0.063% | 674 | 1.46% | 0.059% | 885 | 0.017% | 0.00079% |
| 666 | 0.53% | 0.028% | 802 | 0.048% | 0.0022% | 925 | 0.0043% | 0.00019% |

Table 7: Two initial basis. Left: $tol = 10^{-3}$. Middle: $tol = 10^{-4}$. Right: $tol = 10^{-5}$.

In our next numerical example, the online enrichment is performed for coarse regions that have a cumulative residual that is $\theta$ fraction of the total residual. Assume that the local residuals are

| DOF | $e_a$ | $e_2$ | DOF | $e_a$ | $e_2$ | DOF | $e_a$ | $e_2$ |
|---|---|---|---|---|---|---|---|---|
| 675 | 16.96% | 2.54% | 675 | 16.96% | 2.54% | 900 | 0.54% | 0.023% |
| 863 | 0.63% | 0.027% | 898 | 0.054% | 0.023% | 1087 | 0.011% | 0.00043% |
| 867 | 0.44% | 0.018% | 993 | 0.046% | 0.0015% | 1106 | 0.0050% | 0.00019% |

Table 8: Three initial basis. Left: $tol = 10^{-3}$. Middle: $tol = 10^{-4}$. Right: $tol = 10^{-5}$.

arranged so that
$$r_1 \geq r_2 \geq r_3 \geq \cdots.$$
Then, we only add the basis $\phi_1, \cdots, \phi_k$ for the coarse neighborhoods $\omega_1, \cdots, \omega_k$ such that $k$ is the smallest integer with
$$\theta \sum_{i=1}^{N} r_i^2 \leq \sum_{i=1}^{k} r_i^2.$$

In Table 9, we present numerical results for the last 4 iterations when using 1, 2 and 3 initial basis functions with the tolerance $10^{-4}$ and $\theta = 0.7$. We observe that one can reduce the total number of basis functions compared to the previous case to achieve a similar error. Our conclusions regarding the importance of ONERP condition for $V_{\text{ms}}$ is the same as before.

| DOG | $e_a$ | $e_2$ | DOF | $e_a$ | $e_2$ | DOF | $e_a$ | $e_2$ |
|---|---|---|---|---|---|---|---|---|
| 620 | 1.10% | 0.23% | 450 | 26.60% | 6.92% | 675 | 16.96% | 2.54% |
| 709 | 0.49% | 0.082% | 576 | 1.94% | 0.12% | 827 | 1.02% | 0.045% |
| 787 | 0.050% | 0.0024% | 690 | 0.20% | 0.0099% | 957 | 0.091% | 0.0033% |
| 789 | 0.046% | 0.0022% | 744 | 0.051% | 0.0023% | 987 | 0.048% | 0.0017% |

Table 9: The results using cumulative errors with $\theta = 0.7$, $tol = 10^{-4}$. Left: One initial basis. Middle: Two initial basis. Right: Three initial basis.

### 5.3 Inexpensive online basis construction

In this section, we present numerical results by computing online basis functions in a smaller dimensional spaces. In our first set of numerical examples, the online basis functions are computed in the next $N_0$ snapshot vectors that are eigenvectors of the local spectral problem as discussed earlier. We present numerical results for $N_0 = 10, 20, 40, 50$ in Figure 6. As we observe from this figure that using reduced dimensional spaces for the computation of the online basis function works well only for the first iteration. Later iterations are not very effective. Indeed, because in later iterations, we only add 1 extra snapshot vector to compute the next online basis function. To remedy this, we compute each next iteration of the online basis function by adding more "next" snapshot vectors. In Figure 7, for each next iterate for computing online basis functions, we consider 40 next snapshots and use neighboring regions for $\omega_i$. In this case, the convergence is fast and its behavior is comparable to the use of the whole local snapshot space for the computation of online basis functions.

We would like to remark that the online basis function has no sparsity pattern in the snapshot space and, thus, we could not apply sparsity techniques.

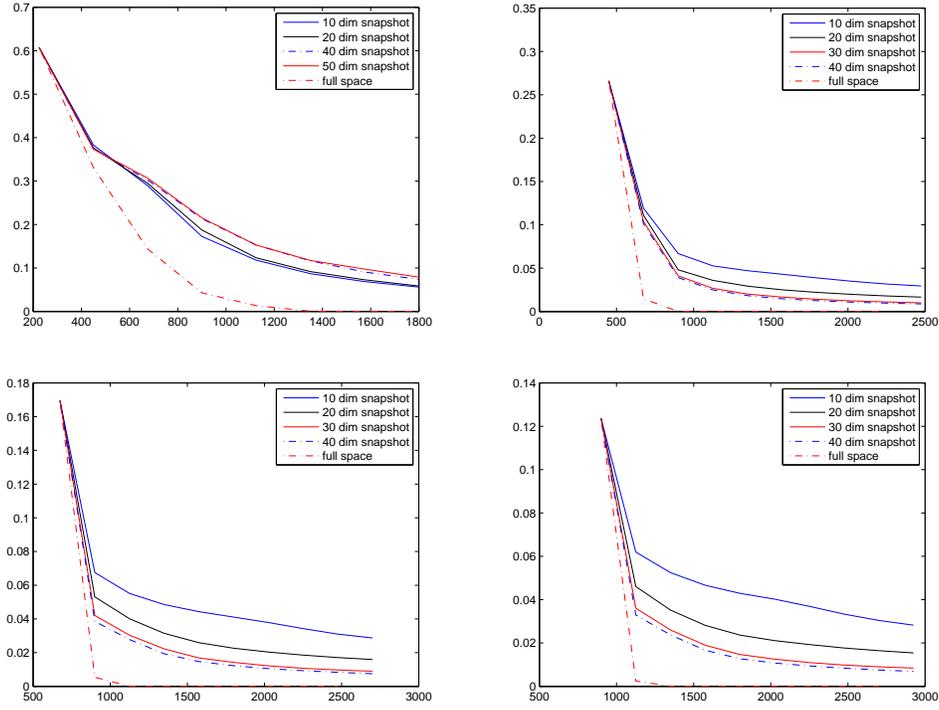

Figure 6: Top-Left: One basis. Top-Right: Two basis. Bottom-Left: Three basis. Bottom-Right: Four basis. Along $x$-axis: Dimensions of $V_{\text{ms}}$. Along $y$-axis: Relative energy errors.

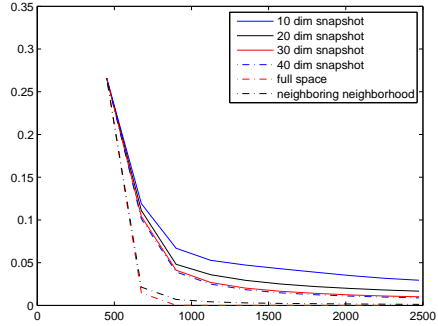

Figure 7: Using the basis for the neighbouring neighborhoods. Along $x$-axis: Dimensions of $V_{\text{ms}}$. Along $y$-axis: Relative energy errors.

**Remark 5.1.** *In all above examples, we consider a snapshot space that consists of local eigenvectors. It can be that in different regions, one needs to use different eigenvalue problems. For example, if the heterogeneities are very localized, one can use polynomial basis functions away from heterogeneous regions. In this case, we consider an adaptive strategy that can identify which class of basis functions to use in a given region. We have considered several basis sets and used the residual to decide which set to use in a given $\omega_i$ for the permeability field shown in Figure 4 with $16 \times 16$ coarse mesh and* each coarse block is subdivided into $16 \times 16$ fine blocks. The set is identified using a few initial basis functions (e.g., using multiscale basis functions). Then, for each set, we consider the residual using only a few basis and estimate the error. We choose the set that gives the largest reduction in the residual. Our numerical results show that by selecting appropriate class of basis functions in each region, we can improve the accuracy of GMsFEM.

# 6 Conclusions

In this paper, we consider a residual-based multiscale basis construction within GMsFEM. The main idea of the proposed method is to construct the online basis functions by solving local problems based on a computed residual. In particular, in each coarse region, an online multiscale basis function is constructed by solving local problems with a right hand side that is a residual computed at the current solution iterate. We show that the offline space needs to satisfy ONERP condition in order to guarantee that adding online basis function will decrease the error independent of the contrast and small scales. The online basis functions account for global effects that are missing in local GMsFEM basis functions. However, the first several GMsFEM basis functions are needed in order to guarantee that the online basis functions will decrease the error independent of the contrast. This method is applied in conjunction with an adaptivity where online basis functions are added in selective regions. The overall procedure results to a local multiscale approach where one can adaptively select regions and compute multiscale basis functions without resorting to global solves. We test our approaches on several examples and present some representative numerical results. Our numerical results show that with the offline spaces that satisfy ONERP, one can achieve a rapid decay of the error. We propose some strategies to reduce the computational cost associated with calculating the online basis functions.